\documentclass[a4paper,11pt]{article}

\usepackage{MyLaTeX,MyArticle}
\input xy
\xyoption{all}

\title{Simple Modules for Groups with Abelian Sylow $2$-Subgroups are Algebraic}
\author{David A.~Craven, University of Oxford}
\date{December 2007}

\begin{document}
\maketitle

\noindent An \emph{algebraic} module is a $KG$-module that satisfies a polynomial with integer coefficients, with addition and multiplication given by direct sum and tensor product. In this article we prove that if $G$ is a group with abelian Sylow $2$-subgroups and $K$ is a field of characteristic $2$, then every simple $KG$-module is algebraic.

\section{Introduction}

The concept of an algebraic module originated with Alperin \cite{alperin1976b}. It can be thought of as an attempt to distinguish those modules whose tensor powers are `nice' from those whose tensor powers are `uncontrollable'. Define a module to be \emph{algebraic} if it satisfies a polynomial with integer coefficients in the Green ring. This is equivalent to the statement that for a module $M$ there is a finite list of indecomposable modules $M_1,\dots,M_n$ such that for every $i$, the module $M^{\otimes i}$ is isomorphic to a sum of the $M_j$.

The main theorem in this article is the following.

\begin{thma}\label{mainthm} Let $G$ be a group with abelian Sylow $2$-subgroups, and let $M$ be a simple $KG$-module, where $K$ is a field of characteristic $2$. Then $M$ is algebraic.\end{thma}

This result is in some sense orthogonal to a recent theorem \cite{cekl2008}, which states that if $B$ is a block with Klein four defect group, then all simple $B$-modules are algebraic. This result implies that there are only finitely many Puig equivalence classes of block with Klein four defect group. Both results can be subsumed under the following conjecture.

\begin{conja}\label{mainconj} Let $B$ be a $2$-block with abelian defect group. Then all simple $B$-modules are algebraic.\end{conja}

This conjecture can not be extended to odd primes in general: the principal block of the sporadic simple group $M_{11}$ contains non-algebraic simple modules, such as the 5-dimensional simple modules. However, it appears as though if a group has a block with abelian defect group then the simple modules in that block tend to be algebraic.

\section{Preliminaries}

The theory of algebraic modules is not well-known, and thus we will recall some results from the literature that we need in order to prove our main theorem. We begin with the following easy lemma.

\begin{lem}[{{\cite[Section II.5]{feit}}}]\label{basicprops} Let $M=M_1\oplus M_2$ be a $KG$-module, and suppose that $H_1\leq G\leq H_2$.
\begin{enumerate}
\item $M$ is algebraic if and only if $M_1$ and $M_2$ are algebraic.
\item The module $M_1\otimes M_2$ is algebraic if $M_1$ and $M_2$ are algebraic.
\item If $M_1$ is algebraic then the modules $M_1\res {H_1}$ and $M_1\ind {H_2}$ are algebraic.
\end{enumerate}
\end{lem}

In particular, this proves the following two corollaries.

\begin{cor}\label{greencorr} Let $K$ be a field and let $G$ be a finite group. Let $M$ be an indecomposable $KG$-module, $S$ be a source of $M$, and $N$ be a Green correspondent of $M$. Then the following are equivalent:
\begin{enumerate}
\item $M$ is algebraic;
\item $N$ is algebraic; and
\item $S$ is algebraic.
\end{enumerate}
\end{cor}

\begin{cor}\label{p'index} Let $G$ be a finite group and let $K$ be a field of characteristic $p$. Let $H$ be a normal subgroup of $G$ whose index is prime to $p$, and let $M$ be a $KG$-module. Then $M$ is algebraic if and only if $M\res H$ is algebraic.
\end{cor}

Lemma \ref{basicprops}(ii), together with the fact that every simple module for a direct product is the tensor product of simple modules for each factors yields the following result.

\begin{lem}\label{directprods} Let $K$ be a field of characteristic $p$. Suppose that $G$ and $H$ are finite groups such that every simple module over $K$ is algebraic. Then every simple $K(G\times H)$-module is algebraic.
\end{lem}

We would like to ignore projective modules when analyzing tensor products of modules. The next result provides this.

\begin{prop}\label{quotbyprojs} Let $\ms I$ be an ideal of algebraic modules in the Green ring $a(KG)$, and let $M$ be a $KG$-module. Then $M$ is algebraic in $a(KG)$ if and only if $M+\ms I$ is algebraic in $a(KG)/\ms I$. In particular, if $\ms P$ denotes the ideal consisting of all projective modules, then a $KG$-module $M$ is algebraic if and only if $M + \mathscr P$ is algebraic.
\end{prop}
\begin{pf} Suppose that $M$ is algebraic. Then $M$ satisfies some polynomial in the Green ring, and therefore its coset in any quotient satisfies this polynomial as well. Conversely, suppose that $M+\ms I$ satisfies some polynomial in the quotient $a(KG)/\ms I$. Thus
\[ \sum \alpha_i (M+\ms I)^i=\ms I.\]
This implies that, since $(M+\ms I)^i=M^{\otimes i}+\ms I$, we have
\[ \sum \alpha_i M^{\otimes i} \in \ms I,\]
which consists solely of algebraic modules. Hence there is some polynomial involving only $M$ witnessing the algebraicity of $M$.
\end{pf}

In fact, since all modules with cyclic vertex are algebraic (as there are only finitely many such indecomposable modules for any group) the collection of all modules whose indecomposable summands all have either cyclic or trivial vertex is an ideal of algebraic modules.

We finally need the soluble version of Theorem \ref{mainthm}.

\begin{thm}[Berger \cite{berger1976}]\label{solalg} Let $G$ be a soluble group, and let $M$ be a simple $KG$-module, where $K$ is a field. Then $M$ is algebraic.
\end{thm}

This establishes all of the basic properties of algebraic modules that we need. Our next results are on Clifford theory \cite{clifford1937}. Let $K$ be a field of characteristic $p$, and let $G$ be a finite group with a normal subgroup $H$. If $M$ is a simple $KG$-module, then $M\res H$ is semisimple; let $N$ be one of its summands. The inertia subgroup $L$ of $N$ is the set of all $g\in G$ such that the conjugate module $N^g$ is isomorphic with $N$.

There are simple projective (i.e., homomorphisms $\phi:G\to \PGL_n(K)$) representations $V$ and $W$ of $L$ such that $V\res H=N$ and $W$ is a projective $G/H$-representation, such that $(V\otimes W)\ind G=M$. A projective representation of a finite group $G$ can be thought of as an ordinary representation of a group $\hat G$, which is a central extension of $G$ by a group of $p'$-order. We will consider this viewpoint when dealing with projective representations.

Further preliminaries involve blocks of defect group $V_4$, which are necessary to make dealing with certain groups considerably easier. The Green correspondence for $V_4$ blocks was determined by Erdmann in \cite{erdmann1982} up to a parameter. The following result is an easy corollary of the main theorems of \cite{erdmann1982}.

\begin{lem}\label{v4blocks} Let $B$ be a $2$-block of a finite group $G$, and suppose that the defect group of $B$ is $V_4$. Suppose further that $B$ is real (which is equivalent to at least one of the ordinary characters of $B$ being real). Then a simple $B$-module either has trivial source or is periodic with 2-dimensional source.\end{lem}

We end this section with the main tool used in the proof of Theorem \ref{mainthm}: Walter's characterization of groups with abelian Sylow $2$-subgroups.

\begin{thm}[Walter \cite{walter1969}, Bender \cite{bender1970}]\label{walterthm} Let $G$ be a finite group with abelian Sylow $2$-subgroups. Then there are normal subgroups $H\leq L$ such that
\begin{enumerate}
\item $|H|$ is odd,
\item $|G:L|$ is odd, and
\item $L/H$ is a direct product of simple groups with abelian Sylow $2$-subgroups and an abelian $2$-group.
\end{enumerate} The simple groups with abelian Sylow $2$-subgroups are $\SL_2(2^n)$ for $n\geq 3$, $\PSL_2(q)$ for $q\geq 5$ and $q\equiv 3,5\mod 8$, the Ree groups $^2G_2(q)$ for $q\geq 27$ an odd power of 3, and the Janko group $J_1$.
\end{thm}

\section{Simple Groups}

Walter's characterization of groups with abelian Sylow $2$-subgroups reduces us to proving the result for simple groups and using an inductive argument, based on Clifford theory. However, in order to use Clifford theory, we need to prove Theorem \ref{mainthm} for all odd central extensions of simple groups with abelian Sylow $2$-subgroups first. Luckily, this will not be a problem, since the groups we are interested in have Schur multiplier at most 2. We deal with each class of groups in turn. Let $K$ denote a field of characteristic 2.

We begin with $\SL_2(2^n)$. In \cite{alperin1979}, Alperin determines the tensor products of simple modules in characteristic 2 for the groups $\SL_2(2^n)$, and in particular proves that all simple modules are algebraic. Since the Schur multipler for $\SL_2(2^n)$ is trivial (except for $\SL_2(4)$, for which it is 2), this deals with this class of group.

The groups $\PSL_2(q)$ where $q\equiv 3,5\mod 8$ are studied in \cite{erdmann1977}, where it is proved that if $q\equiv 3\mod 8$ then the sources of the simple modules in the (unique) block of non-cyclic defect group are all trivial. Thus when $q\equiv 3\mod 8$ all simple modules for $\PSL_2(q)$ are algebraic. (The Schur multiplier for these groups is $2$.) When $q\equiv 5\mod 8$, the two non-trivial simple modules in the principal block are periodic. Their source is 2-dimensional, and so can be thought of as a restriction of the natural module $\GL_2(2^n)$ for some $n$. (In fact, in \cite{erdmann1977}, it is shown that $n$ can be taken to be $2$.) Since all simple modules for $\GL_2(2^n)$ are algebraic (by the previous paragraph and Corollary \ref{p'index}), we see that the 2-dimensional sources are likewise algebraic. Thus all simple modules are algebraic when $q\equiv 5\mod 8$.

The groups $^2G_2(q)$ are more complicated. The isomorphism type of the normalizer of a Sylow $2$-subgroup (which is isomorphic to $(C_2\times C_2\times C_2)\rtimes T$, where $T$ is the non-abelian group of order 21) does not depend on the particular $q$. In \cite{landrockmichler1980}, Landrock and Michler proved that the Green correspondents for the simple modules in the principal block of $^2G_2(q)$ are independent of $q$. In particular, this correspondence holds for $^2G_2(3)=\SL_2(8)\rtimes C_3$. Then Corollaries \ref{greencorr} and \ref{p'index}, and the fact that all simple modules for $\SL_2(2^n)$ are algebraic, prove that all simple modules in the principal block are algebraic.

We now use results of Ward in \cite{ward1966}. The remaining simple modules lie in blocks of smaller defect. All blocks with trivial or cyclic defect group are populated solely by algebraic modules, and so the discussion rests with blocks of defect group $V_4$. In \cite[II.7]{ward1966}, Ward notes that the ordinary characters in $2$-blocks of defect 2 are real. By Lemma \ref{v4blocks}, this implies that all simple modules from blocks of defect group $V_4$ are algebraic. (One may also use the main result of \cite{cekl2008}, as mentioned before. That result uses the same method of proof for this group.)

This leaves the simple group $J_1$. In \cite{landrockmichler1978}, the block structure of $J_1$ is given. It has one block of full defect, together with a cyclic block and five projective simple modules. Hence all simple modules are algebraic if and only if all simple modules from the principal block are algebraic. The normalizer of a Sylow $2$-subgroup of $J_1$ is the same as that of $^2G_2(q)$, and the Green correspondent of the 20-dimensional simple module in the principal block of $J_1$ is the same as the Green correspondent of the 12-dimensional simple module for $\SL_2(8)\rtimes C_3$. This leaves the two (dual) 56-dimensional simple modules and the 76-dimensional simple module.

To analyze the algebraicity of these modules requires the use of a computer. The author knows of no way to compute the resulting tensor products by hand, since the dimensions and numbers of indecomposable modules involved make most decompositions highly difficult.

The decompositions of tensor powers of the Green correspondents of the 76-dimensional simple module, and of the sources of the 56-dimensional simple modules, were computed using the computer algebra package Magma, and the results, together with a fwe remarks on the methods needed to compute the decompositions, are given in Section \ref{j1}.

We have the following theorem.

\begin{thm}\label{simplemainthm} Let $G$ be a finite simple group with abelian Sylow $2$-subgroups, and let $K$ be a field of characteristic $2$. Then every simple projective $KG$-representation is algebraic.
\end{thm}

\section{The Induction}

Recall the notation $\Orth_{2'}(G)$ and $\Orth^{2'}(G)$ for the largest normal subgroup of $G$ of odd order and the smallest normal subgroup of odd index respectively.

By Corollary \ref{p'index}, if Theorem \ref{mainthm} is true for all such finite groups with $G=\Orth^{2'}(G)$ then it is true for all finite groups. We actually prove something equivalent but apparently slightly stronger.

\begin{thm} Let $G$ be a finite group with abelian Sylow $2$-subgroups, and let $K$ be a field of characteristic $2$. Let $M$ be a simple projective representation of $G$. Then $M$ is algebraic.
\end{thm}
\begin{pf} Write $H=\Orth_{2'}(G)$, and proceed by induction on $|G/H|$. Let $M$ be a simple projective $KG$-representation. By replacing $G$ with an odd central extension of $G$ (and noting that $G/\Orth_{2'}(G)$ does not change) we may assume that $M$ is a simple $KG$-module. Let $N$ be a simple summand of $M\res H$ and let $L$ denote the inertia subgroup of $N$. If $L=G$, then write $V$ for the projective $KG$-representation such that $V\res H=N$. If $P$ is a Sylow $2$-subgroup of $G$, then $V\res{HP}$ is a simple projective $K(HP)$-representation, which is therefore algebraic, since $HP$ is soluble (by Theorem \ref{solalg}). Thus $V$ is algebraic.

We can also see that $W$ is algebraic, since it is a projective $G/H$-representation; by Theorem \ref{walterthm}, $G/H$ is a direct product of a $2$-group $P$ and simple groups with abelian Sylow $2$-subgroups. A central extension of $G/H$ of odd order must be a group $A\times B$, where $A$ is a direct product of simple groups and $B$ is an odd central extension of $P$. By Theorem \ref{simplemainthm}, all simple modules for each factor of the product $A$ are algebraic, and by Theorem \ref{solalg}, all simple modules for $B$ are algebraic. Finally, by Lemma \ref{directprods}, all simple modules for all odd central extensions of $G/H$ are algebraic, as required.

Now suppose that $L<G$. Then there is a simple $KL$-module $M'$ such that $M'\ind G=M$. The subgroup $L$ is a group with abelian Sylow $2$-subgroups, and since $H\leq L<G$, we must have $|L/\Orth_{2'}(L)|<|G/\Orth_{2'}(G)|$. By induction all projective simple $KL$-representations are algebraic. Thus $M'$ is algebraic, and so is $M$, as required.\end{pf}

\section{Computing with $J_1$}
\label{j1}
Write $G$ for the normalizer of a Sylow $2$-subgroup $P$ of $J_1$, which has order 168, and let $K$ be a field of characteristic 2. The Green correspondents of the simple modules for $J_1$ are given in \cite{landrockmichler1978}. There are therefore two indecomposable modules for $G$ that need to be analyzed. Unfortunately, it does not seem possible at this juncture to provide theoretical reasons why these two indecomposable modules are algebraic, and to manually compute the tensor products is a daunting task unless one is aided by a computer.

It is easy to see that there are five simple $KG$-modules if $K$ contains a cube root of unity, and four simple modules otherwise: the trivial module, two other dual 1-dimensional modules which amalgamate over $\GF(2)$, and two dual 3-dimensional modules.

The first module to be considered is the 12-dimensional Green correspondent $A$ of the 76-dimensional simple module of $J_1$ in the principal block. Its Green correspondent does not require a cube root of unity to be realized, and so we let $K$ be the field $\GF(2)$. Therefore $G$ has four simple modules: the trivial module $K$, the 2-dimensional simple (but not absolutely simple) module $W_1$, and the two dual 3-dimensional simple modules $W_2$ and $W_2^*$. Finally, if $M$ is a module, write $\proj(M)$ for the projective cover of $M$, and write $n\cdot M$ for the $n$-fold direct sum of $M$ with itself.

The tensor square of $A$ is given by
\[ A\otimes A=2\cdot C_1\oplus C_2,\]
where $C_1$ is a 28-dimensional module with $\Omega(C_1)=C_1$ and $C_2$ is a non-periodic 88-dimensional indecomposable module. The decomposition of $A\otimes C_1$ is given by
\[ A\otimes C_1=2\cdot \proj(K)\oplus 2\cdot \proj(W_1)\oplus 4\cdot\proj(W_2)\oplus 4\cdot \proj(W_2^*)\oplus 2\cdot C_3,\]
where $C_3$ is a 56-dimensional module with $\Omega(C_3)=C_3$. Finally in this direction, we have the decomposition
\[ A\otimes C_3=2\cdot \proj(K)\oplus 3\cdot \proj(W_1)\oplus 8\cdot\proj(W_2)\oplus 8\cdot \proj(W_2^*)\oplus 4\cdot C_1\oplus 2\cdot C_3.\]

Moving on to the product of $A$ and $C_2$, we have
\[ A\otimes C_2=4\cdot \proj(K)\oplus 3\cdot A_1\oplus 4\cdot\proj(W_1)\oplus 14\cdot \proj(W_2)\oplus 14\cdot \proj(W_2^*)\oplus 2\cdot C_1\oplus 2\cdot C_3\oplus C_4\oplus C_5,\]
where $C_4$ is a 12-dimensional indecomposable module and where $C_5$ is a periodic 84-dimensional indecomposable module. Decomposing more tensor products, we have
\[ A\otimes C_4=\proj(K)\oplus \proj(W_1)\oplus 2\cdot\proj(W_2)\oplus 2\cdot \proj(W_2^*)\oplus C_6,\]
where $C_6$ is a 24-dimensional indecomposable module. The next decomposition is
\[ A\otimes C_6=2\cdot C_4\oplus 2\cdot \proj(K)\oplus 2\cdot \proj(W_1)\oplus 4\cdot \proj(W_2)\oplus 4\cdot \proj(W_2^*)\oplus C_7,\]
where $C_7$ is another 24-dimensional indecomposable module. Decomposing the next tensor product gives
\[ A\otimes C_7=2\cdot \proj(K)\oplus 2\cdot \proj(W_1)\oplus 4\cdot \proj(W_2)\oplus 4\cdot \proj(W_2^*)\oplus C_8,\]
where $C_8$ is a 48-dimensional indecomposable module. Finally, we have
\[ A\otimes C_8=2\cdot C_4\oplus 2\cdot C_7\oplus 3\cdot \proj(K)\oplus 3\cdot \proj(W_1)\oplus 9\cdot\proj(W_2)\oplus 9\cdot \proj(W_2^*).\]

Thus the remaining module to deal with is $C_5$. This decomposes as
\[ A\otimes C_5=5\cdot \proj(K)\oplus 5\cdot \proj(W_1)\oplus 15\cdot \proj(W_2)\oplus 15\cdot \proj(W_2^*)\oplus C_9,\]
where $C_9$ is a 168-dimensional periodic module, and this is followed by
\[ A\otimes C_9=10\cdot \proj(K)\oplus 10\cdot \proj(W_1)\oplus 30\cdot \proj(W_2)\oplus 30\cdot \proj(W_2^*)\oplus 2\cdot C_5\oplus C_{10},\]
where $C_{10}$ is another 168-dimensional periodic module. Next,
\[ A\otimes C_{10}=10\cdot \proj(K)\oplus 10\cdot \proj(W_1)\oplus 30\cdot \proj(W_2)\oplus 30\cdot \proj(W_2^*)\oplus C_{11},\]
where $C_{11}$ is a 336-dimensional periodic module, and finally
\[ A\otimes C_{11}=21\cdot \proj(K)\oplus 21\cdot \proj(W_1)\oplus 63\cdot \proj(W_2)\oplus 63\cdot \proj(W_2^*)\oplus 2\cdot C_{5}\oplus 2\cdot C_{10},\]
proving that $A$ is algebraic. Therefore by Corollary \ref{greencorr}, the corresponding 76-dimensional simple module for $J_1$ is algebraic.

The last two modules to analyze are the two dual 56-dimensional simple modules for $J_1$. In this case, we really do need to extend our field to $K=\GF(4)$, and we do so. Let $B$ denote the source of one of these modules, an 8-dimensional non-periodic $KP$-module. The module $B^{\otimes 2}$ is indecomposable, but
\[ B^{\otimes 3}=2\cdot D_1\oplus \bigoplus_{i=1}^3 D_{2,i}\oplus Y\oplus 48\cdot\proj(K),\]
where $D_1$ is a periodic 8-dimensional module, the $D_{2,i}$ are periodic 28-dimensional modules, and $Y$ is a sum of 4-dimensional permutation modules with cyclic vertex (which can be ignored by Proposition \ref{quotbyprojs}). The tensor product of $B$ and $D_1$ is simply a sum of modules with cyclic or trivial vertex, and so we consider the tensor product of $B$ and $D_{2,i}$. In this case,
\[ B\otimes D_{2,i}=D_{3,i}\oplus 21\cdot \proj(K),\]
where the $D_{3,i}$ are (pairwise non-isomorphic) 56-dimensional indecomposable modules. Next,
\[ B\otimes D_{3,i}=2\cdot D_{2,i}\oplus D_{4,i}\oplus 42\cdot \proj(K),\]
where the $D_{4,i}$ are (pairwise different) 56-dimensional indecomposable modules. The module $D_{4,i}\otimes B$ behaves similarly to $D_{2,i}\otimes B$, and indeed
\[ B\otimes D_{4,i}=D_{5,i}\oplus 42\cdot \proj(K),\]
where the $D_{5,i}$ are 112-dimensional indecomposable modules, and
\[ B\otimes D_{5,i}=2\cdot D_{3,i}\oplus 2\cdot D_{5,i}\oplus 91\cdot \proj(K).\]
This implies that $B$ is algebraic, since we have decomposed all possible tensor products.

This confirms that all simple modules in the principal block of $J_1$ are algebraic, providing the last piece in the proof of Theorem \ref{mainthm}.
\\ \qquad

We will end with a brief discussion on how to produce the decompositions listed above. The standard method in Magma is to use the in-built command \texttt{IndecomposableSummands} which, as its name implies, returns the indecomposable summands of a given module. However, this method quickly runs out of memory when dealing with modules of large dimension. (It appears to be dependent rather on the number of submodules of a module, which makes this method particularly bad for $p$-groups.) For example, the 896-dimensional module $D_{5,1}\otimes B$ over $\GF(4)$ would be difficult to decompose using this command.

Instead, we make the observation that the majority of the summands in tensor products of indecomposable modules for $p$-groups are projective. Thus we need only find \emph{submodules} of a module that are projective, and we immediately know that they split off. This removal of projective summands can be easily done for $p$-groups with very little memory requirements, using the fact that if a module for a $p$-group has a free summand, then it is likely that the submodule generated by a random element is free.

This method will not work very well (indeed, at all) for a group that is not a $p$-group. For example, we are interested in the tensor powers of the Green correspondent, and so we need another method. It is possible to construct the set of all module homomorphisms between a given module and a projective $KG$-module. Choosing random elements from the set of homomorphisms yields a (probabilistic) method of removing a projective summand. This requires the construction of all projective $KG$-modules.

\bibliography{references}

\end{document}